\documentclass[11pt]{article}

\usepackage{amscd}
\usepackage{amsfonts}
\usepackage{amsmath}
\usepackage{amssymb}
\usepackage{amsthm}
\usepackage{latexsym}
\usepackage{mathrsfs}
\usepackage[mathcal]{euscript}
\usepackage{xcolor}
\usepackage{array}
\usepackage{multirow}

\font\smallit=cmti10
\font\smalltt=cmtt10

\newtheorem{thm}{Theorem}
\newtheorem{prop}[thm]{Proposition}

\newtheorem{cor}[thm]{Corollary}
\newtheorem{rem}[thm]{Remark}

\newtheorem{defn}[thm]{Definition}
\newenvironment{Pf}{\noindent{\bf Proof. }}{\hfill $\blacksquare$ \\}

\def\ds{\displaystyle}
\def\ts{\textstyle}
\def\lc{\lceil}
\def\rc{\rceil}

\def\mod{\!\!\pmod}
\def\imp{\rightarrow}
\def\ov{\overline}
\def\pr{\prime}
\def\sm{\setminus}

\def\l{\lambda}

\def\s{\sigma}

\def\D{\Delta}

\def\L{\mathcal L}
\def\N{\mathbb N}
\def\Rad{\tt Rad}
\def\S{\mathrm S}
\def\T{\mathrm T}
\def\Z{\mathbb Z}

\begin{document}
\title{\bf On the two-colour Rado number for $\sum_{i=1}^m a_ix_i=c$}
\author{
{\bf Ishan Arora}\thanks{{\smallit CEO, Cerebro Limited, Grey Lynn, Auckland, 1021 New Zealand} \newline {\smalltt e-mail:\,ishanarora@gmail.com}}\:\:\thanks{Work done while at Department of Mathematics, Indian Institute of Technology Delhi} 
\qquad 
{\bf Srashti Dwivedi}\thanks{{\smallit Department of Mathematics, Alliance University, Anekal, Bengaluru -- 562106, India} \newline {\smalltt e-mail:\,srashti.dwivedi@alliance.edu.in}}\:\:\thanks{Work done while at Department of Mathematics, Indian Institute of Technology Delhi} 
\qquad 
{\bf Amitabha Tripathi}\thanks{{\smallit Department of Mathematics, Indian Institute of Technology, Hauz Khas, New Delhi -- 110016, India} \newline {\smalltt e-mail:\,atripath@maths.iitd.ac.in}}\:\:\thanks{Corresponding author} 
}
\date{}
\maketitle

\begin{abstract}
\noindent Let $a_1,\ldots,a_m$ be nonzero integers, $c \in \Z$ and $r \ge 2$. The {\tt Rado number} for the equation
\[ \sum_{i=1}^m a_ix_i = c \]
in $r$ colours is the least positive integer $N$ such that any $r$-colouring of the integers in the interval $[1,N]$ admits a monochromatic solution to the given equation. We introduce the concept of $t$-distributability of sets of positive integers, and determine exact values whenever possible, and upper and lower bounds otherwise, for the Rado numbers when the set $\{a_1,\ldots,a_{m-1}\}$ is $2$-distributable or $3$-distributable, $a_m=-1$, and $r=2$. This generalizes previous works by several authors. 
\end{abstract}
\vskip 10pt

\noindent {\bf Keywords:} $r$-colouring, monochromatic solution, valid colouring, Rado number, $r$-distributability \\
{\bf 2010 Mathematics Subject Classification:} 05C55 (primary), 05D10 (secondary) \\

\section{Introduction}

I. Schur \cite{Sch16} discovered his celebrated result that bears his name in 1916 while attempting to prove Fermat's ``last" theorem. Schur's theorem states that for every positive integer $r$, there exists a positive integer $n=n(r)$ such that for every $r$-colouring of the integers in the interval $[1,n]$, there exists a monochromatic solution to the equation $x+y=z$. In other words, for each positive integer $r$, and for every function $\chi: \{1,\ldots,n\} \imp \{1,\ldots,r\}$, there exist integers $x, y \in \{1,\ldots,n\}$ with $x+y \in \{1,\ldots,n\}$ such that ${\chi}(x)=\chi(y)=\chi(x+y)$. The least such $n=n(r)$ is denoted by ${\tt s}(r)$ in honour of Schur, and these are called the Schur numbers. The only exact values of ${\tt s}(r)$ known are ${\tt s}(1)=2$, ${\tt s}(2)=5$, ${\tt s}(3)=14$, ${\tt s}(4)=45$, and ${\tt s}(5)=160$. 

Schur's theorem was generalized in a series of results in the 1930's by Rado \cite{Rad33a,Rad33b,Rad45} leading to a complete resolution to the following problem: characterize systems of linear homogeneous equations with integral coefficients $\L$ such that for a given positive integer $r$, there exists a least positive integer $n=\Rad(\L;r)$ for which every $r$-colouring of the integers in the interval $[1,n]$ yields a monochromatic solution to the system $\L$. There has been a growing interest in the determination of the Rado numbers $\Rad(\L;r)$, particularly when $\L$ is a single equation and $r=2$. 
 
Beutelspacher \&  Brestovansky \cite{BB82} proved that the $2$-colour Rado number for the equation $x_1+\cdots+x_{m-1}-x_m=0$ equals $m^2-m-1$ for each $m \ge 3$. The $2$-colour Rado numbers for the equation $a_1x_1+a_2x_2+a_3x_3=0$ has been completely resolved. In view of Rado's result on single linear homogeneous regular equations, the only equations to consider are when $a_1+a_2=0$ and when $a_1+a_2+a_3=0$. The first case was completely resolved by Burr \& Loo \cite{BL}, and the second case by Gupta, Thulasi Rangan \& Tripathi \cite{GTT15}. Jones \& Schaal \cite{JS01} determined the $2$-colour Rado number for $a_1x_1+\cdots+a_{m-1}x_{m-1}-x_m=0$ when $a_1,\ldots,a_{m-1}$ are positive integers with $\min\{a_1,\ldots,a_{m-1}\}=1$. Hopkins \& Schaal \cite{HS05} resolved the problem in the case $\min\{a_1,\ldots,a_{m-1}\}=2$ and gave bounds in the general case which they conjectured to hold; this was proved by Guo \& Sun in \cite{GS08}. They showed that the $2$-colour Rado number for the general case equals $as^2+s-a$, where $a=\min\{a_1,\ldots,a_{m-1}\}$ and $s=a_1+\cdots+a_{m-1}$. Kosek \& Schaal \cite{KS01} determined the $2$-colour Rado number for the non-homogeneous equation $x_1+\cdots+x_{m-1}-x_m=c$ when $c<m-2$ and when $c>(m-1)(m-2)$. In all other cases, they have given at least upper and lower bounds. 

Let $a_1,\ldots,a_m$ be nonzero integers, $c \in \Z$ and $r \in \N$. The {\tt Rado number} for the equation
\begin{equation} \label{geneqn}
\sum_{i=1}^m a_ix_i = c
\end{equation}
in $r$ colours, denoted by ${\tt Rad}_r \big(a_1,\ldots,a_m;c\big)$, is the least positive integer $N$ such that any $r$-colouring of the integers in the interval $[1,N]$ admits a monochromatic solution to Eq.(\ref{geneqn}). This means that for {\it any\/} $\chi: [1,N] \imp \{1,\ldots,r\}$ there exists $x_1,\ldots,x_m \in [1,N]$ satisfying Eq.(\ref{geneqn}) and such that ${\chi}(x_1)=\cdots={\chi}(x_m)$. 

The case $r=1$ is trivial. For the most of the remainder of this paper, we deal with the case $r=2$. It is convenient to denote ${\tt Rad}_2 \big(a_1,\ldots,a_m;c\big)$ more simply by ${\tt Rad} \big(a_1,\ldots,a_m;c\big)$. In this paper, we provide exact results whenever possible and upper and lower bounds otherwise for these Rado numbers, mostly in the case $r=2$. These results cover all integral values of $c$ but hold for a restricted collection of coefficients $\{a_1,\ldots,a_{m-1}\}$; we always assume $a_m=-1$. We call this restricted family ``distributed", and introduce the basic concepts of ``$r$-distributability" in Section 2. We give the main results in Section 3. A summary of our results is provided in Table 1.  

\begin{table}[ph]
\centering
\renewcommand{\arraystretch}{1.2}
\begin{tabular}{|c|c|c|c|} \hline
range for $c$ & $r$ & restriction on the set & Rado number \\ \hline
\multirow{3}{*}{$(-\infty,\S-1)$} & \multirow{3}{*}{2} & no restriction         & $\ge (\S-1-c)(\S+2)+1$  \\ 
                                                          &                                   & \multirow{2}{*}{$2$-distributable} & $(\S-1-c)(\S+2)+1$ \\
                                                          &                                   &                                & (Theorem \ref{c<S-1}) \\ \hline
\multirow{2}{*}{$\S-1$} & \multirow{2}{*}{all} & \multirow{2}{*}{no restriction} & 1 \\
                                            &                                     &                         & (Theorem \ref{c<S-1}) \\ \hline
\multirow{2}{*}{$[\S,2\S-3]$} & \multirow{2}{*}{2} & \multirow{2}{*}{$2$-distributable} & $\le \S+1$ \\
                                                      &                                    &                                                               & (Theorem \ref{[S,2S-1]}) \\ \hline
\multirow{2}{*}{$2\S-2$} & \multirow{2}{*}{all} & \multirow{2}{*}{no restriction} & 2 \\
                                              &                                     &                                                         & (Theorem \ref{[S,2S-1]}) \\ \hline
\multirow{2}{*}{$2\S-1$} & \multirow{2}{*}{2} & \multirow{2}{*}{$2$-distributable} & 3 \\
                                              &                                   &                                                                & (Theorem \ref{[S,2S-1]}) \\ \hline
$\big[(\l-1)\S,\l\S-\l\big)$ & \multirow{4}{*}{2} & no restriction        & $\ge \l+\mu$ \\
$(3 \le \l \le \S)$                  &                                   & $2$-distributable & $\l+1$ if $\mu=1$ \\                                                                
$c=\l(\S-1)-\mu$                 &                                   & \multirow{2}{*}{$3$-distributable} & $\l+\mu$ if $2 \le \mu \le \S-\l$ \\
$(1 \le \mu \le \S-\l)$         &                                   &                                                                & (Theorem \ref{lambda_mu_1}) \\ \hline
\multirow{2}{*}{$\l(\S-1)$} & \multirow{2}{*}{all} & \multirow{2}{*}{no restriction} & $\l$ \\
                                                 &                                     &                                                        & (Theorem \ref{lambda_mu_1}) \\ \hline
$\ds\bigcup_{3 \le \l \le \big\lc\frac{\S+1}{2}\big\rc} \big(\l\S-\l,\l\S\big)$ & \multirow{4}{*}{2} & \multirow{4}{*}{$3$-distributable} 
                                                                                                                                      & \multirow{2}{*}{$\le 2\l-\mu$} \\
$c=\l \S-\mu$ &  &  & (Theorem \ref{lambda_mu_2}) \\
$(1 \le \mu \le \l-1)$  &  &  &  \\ \hline
\multirow{4}{*}{$\big(\S-1,\S(\S-1)\big)$} & \multirow{4}{*}{2} & \multirow{2}{*}{no restriction} & $\ge \left\lc\frac{1+c(\S+2)}{{\S}^2+\S-1}\right\rc$ if $c>\S-1$  \\ 
& & & (Theorem \ref{>S(S-1)}) \\ 
& & \multirow{2}{*}{$3$-distributable} & $\le \S+1$ \\ 
& & & (Theorem \ref{S+1}) \\ \hline
\multirow{2}{*}{$\big(\S(\S-1),\infty\big)$} & \multirow{2}{*}{2} & \multirow{2}{*}{$2$-distributable} & $\left\lc\frac{1+c(\S+2)}{{\S}^2+\S-1}\right\rc$ if $c>\S(\S-1)$ \\ 
& & & (Theorem \ref{>S(S-1)}) \\ \hline
\end{tabular}
\caption{\bf Summary of results on ${\tt Rad}_r \big(a_1,\ldots,a_{m-1},-1;c\big)$, $c\,\S$ even, $S=\sum_{i=1}^{m-1} a_i$}
\end{table}

\section{Distributable Sets}

The results of our paper hold whenever the collection of positive integers $a_1,\ldots,a_{m-1}$ is distributable. In this section, we introduce the notion of $t$-distributability, for each positive integer $t$. The case $t=2$ corresponds to finite ``complete" sequences, which are well known. 

\begin{defn}
A sequence $\{u_n\}_{n \ge 1}$ of positive integers is said to be {\tt complete} if every positive integer can be expressed as a sum of distinct terms of some finite subsequence $\{u_{n_k}\}_{k \ge 1}$.  
\end{defn}

\begin{prop} {\bf (Brown, 1961)} \label{complete} \\
A sequence $\{u_n\}_{n \ge 1}$ of positive integers is {\tt complete} if and only if 
\begin{equation}
u_1 = 1, \quad u_n \le s_{n-1}+1 \;\;\text{for}\;\; n \ge 2, 
\end{equation}
where $s_k=u_1+\cdots+u_k$ for $k \ge 1$. 
\end{prop}

\begin{defn}
Let $\{a_1,\ldots,a_k\}$ be a multiset of positive integers and $\{{\s}_1,\ldots,{\s}_t\}$ be a multiset of nonnegative integers such that $\sum_{i=1}^t {\s}_i =\sum_{i=1}^k a_i$. The number of partitions of $\{1,\ldots,k\}$ into sets $S_1,\ldots,S_t$ {\rm (}some possibly empty{\rm )} such that $\sum_{i \in S_j} a_i = {\s}_j$ for $j \in \{1,\ldots,t\}$ is called a {\tt Set Distribution Coefficient}, and denoted by $\left[ \begin{array}{c} a_1,\ldots, a_k \\ {\s}_1,\ldots, {\s}_t \end{array}\right]$. 
\end{defn}

Suppose the sets $S_1,\ldots,S_t$ partition $\{1,\ldots,k\}$ such that $\sum_{i \in S_j} a_i={\s}_j$ for $1 \le j \le t$. Then $k$ belongs to exactly one of $S_1,\ldots,S_t$, say $k \in S_1$ after renumbering. Then $\sum_{i \in S_j} a_i={\s}_j$ for $2 \le j \le t$ and $\sum_{i \in S_1 \sm \{k\}} a_i={\s}_1-a_k$. Thus we have the identity 
\begin{equation} \label{SD_identity}
\left[ \begin{array}{c} a_1,\ldots, a_k \\ {\s}_1,\ldots, {\s}_t \end{array}\right] = \left[ \begin{array}{c} a_1,\ldots, a_{k-1} \\ {\s}_1-a_k,{\s}_2,\ldots, 
{\s}_t  \end{array}\right] + \left[ \begin{array}{c} a_1,\ldots, a_{k-1} \\ {\s}_1,{\s}_2-a_k,\ldots, {\s}_t  \end{array}\right] + \cdots + \left[ \begin{array}{c} a_1,\ldots, a_{k-1} \\ {\s}_1,{\s}_2,\ldots, {\s}_t-a_k \end{array}\right].
\end{equation}
 
\begin{defn}
We say a multiset $\{a_1,\ldots,a_k\}$ of positive integers is {\tt $t$-distributable} if $\left[ \begin{array}{c} a_1,\ldots, a_k \\ {\s}_1,\ldots, {\s}_t \end{array}\right] > 0$ for all choices of multisets $\{{\s}_1,\ldots,{\s}_t\}$ of nonnegative integers for which $\sum_{i=1}^t {\s}_i = \sum_{i=1}^k a_i$.
\end{defn}

\begin{rem}
Any multiset of positive integers is $1$-distributable, and any $t$-distributable multiset is $j$-distributable for $j \in \{1,\ldots,t-1\}$. 
\end{rem}

\begin{prop}
Let $A=\{a_1,\ldots,a_k\}$ be a multiset of positive integers. With $a_i \le a_{i+1}$ for $i \in \{1,\ldots,k-1\}$, $A$ is $t$-distributable if and only if 
\begin{equation} \label{t-distributable}
a_i \le \left\lc\frac{s_i}{t}\right\rc \;\;\text{for}\;\; i \in \{1,\ldots,k\} 
\end{equation} 
where $s_i=a_1+a_2+\cdots+a_i$ for $i \in \{1,\ldots,k\}$. 
\end{prop}

\begin{Pf}
Suppose $A=\{a_1,\ldots,a_k\}$ is $t$-distributable, with $a_1 \le a_2 \le \cdots \le a_k$, and let $\sum_{i=1}^k a_i=\s$. Therefore, for each $i \in \{1,\ldots, k\}$, 
\[ \left[ \begin{array}{c} a_1,\ldots, a_k \\ a_i-1,\ldots, a_i-1, {\s}-(t-1)(a_i-1) \end{array} \right] > 0. \]
Since the terms $a_i,\ldots,a_k$ cannot contribute to the sum $a_i-1$, it follows that 
\begin{equation} \label{t-distributable_ineq}
{\s} - (t-1)(a_i-1) \ge a_i + a_{i+1} + \cdots + a_k. 
\end{equation}
Now Eq.(\ref{t-distributable_ineq}) is equivalent to $a_1+\cdots+a_{i-1} \ge (t-1)(a_i-1)$. Adding $a_i$ to both sides, we have 
\[ s_i \ge ta_i - (t-1) > t(a_i-1). \]
Hence $a_i \le \left\lc\frac{s_i}{t}\right\rc$, proving the necessity of the condition. 

Conversely, we must show that any sequence $\{a_i\}_{i=1}^k$ satisfying the condition in Eq.(\ref{t-distributable}) is $t$-distributable. We induct on $k$. The base case $k=1$ is obvious, and we assume that all sequences satisfying the condition in Eq.(\ref{t-distributable}) is $t$-distributable whenever the sequence has fewer than $k$ terms. 

Consider any nondecreasing sequence $\{a_i\}_{i=1}^k$ of positive integers satisfying the condition in Eq.(\ref{t-distributable}), and let ${\s}_1,\ldots,{\s}_t$ be any nondecreasing sequence of nonnegative integers with sum $\s=\sum_{i=1}^k a_i$. Since the subsequence $\{a_i\}_{i=1}^{k-1}$ also satisfies the condition in  Eq.(\ref{t-distributable}) and ${\s}_1,\ldots,{\s}_{t-1},{\s}_t-a_k$ has sum $\s-a_k=\sum_{i=1}^{k-1} a_i$, by the induction hypothesis we have 
\[ \left[ \begin{array}{c} a_1,\ldots, a_{k-1} \\ {\s}_1,\ldots, {\s}_{t-1}, {\s}_t-a_k \end{array} \right] > 0. \]
From Eq.(\ref{SD_identity}) it follows that 
\[ \left[ \begin{array}{c} a_1,\ldots, a_{k-1},a_k \\ {\s}_1,\ldots, {\s}_{t-1}, {\s}_t \end{array} \right] > 0, \]
proving the sufficiency of the condition. 
\end{Pf}

\begin{rem}
The condition Eq.(\ref{t-distributable}) implies, and is implied by, the condition $a_{t-1}=1$. The last equality forces $a_1=\cdots=a_{t-2}=1$. 
\end{rem}

\section{Main Results}

For any collection of nonzero integers $a_1,\ldots,a_m$, and for $c \in \Z$ and $r>1$, the Rado number ${\tt Rad}_r \big(a_1,\ldots,a_m;c\big)$ is the smallest positive integer $R$ for which every $r$-colouring of $[1,R]$ contains a monochromatic solution to Eq.(\ref{geneqn}). By assigning the colour of $x_i$ in the solution to Eq.(\ref{geneqn}) to $x_i-1$, we note that this is equivalent to determining the smallest positive integer $R$ for which every $r$-colouring of $[0,R-1]$ contains a monochromatic solution to 
\begin{equation} \label{revgeneqn}
\sum_{i=1}^m a_ix_i = c-\s,
\end{equation}
where $\s=\sum_{i=1}^m a_i$. 

\begin{thm} \label{DNE}
Let $a_1,\ldots,a_m$ be nonzero integers, $c \in \Z$ and $r>1$. If $\sum_{i=1}^m a_i=\s$, $\text{lcm}(2,3,\ldots,r)=L$, and $\gcd(\s,L) \nmid c$, then ${\tt Rad}_r (a_1,\ldots,a_m;c)$ does not exist. 
\end{thm}

\begin{Pf}
Let $R \ge 1$. We exhibit an $r$-coloring of $[1,R]$ which contains no monochromatic solution to Eq.(\ref{revgeneqn}). Since $\gcd(\s,L) \nmid c$, there exists a prime $p$ that divides both $\s$ and $L$ but does not divide $c$. Observe that $p \le r$. 

Define the colouring $\chi: \N \imp \{1,\ldots,p\}$ by ${\chi}(i) \equiv i\bmod{p}$.  We show that this colouring does not admit a monochromatic solution to Eq.(\ref{geneqn}). By way of contradiction, assume there exists a monochromatic solution $x_1,\ldots,x_m$. Since ${\chi}(x_i)$ is constant, $x_i \equiv x_j\mod{p}$ for $i \ne j$ by definition of $\chi$. Write $x_i \equiv x_0\mod{p}$. Then $c=\sum_i^m a_ix_i \equiv x_0 \sum_{i=1}^m a_i = \s x_0 \equiv 0\mod{p}$, which is a contradiction to our assumption. 
\end{Pf}

\begin{cor} \label{dne}
Let $a_1,\ldots,a_m$ be nonzero integers, $c \in \Z$ and $r=2$. If $\sum_{i=1}^m a_i=\s$ is even and $c$ is odd, then ${\tt Rad}_2 (a_1,\ldots,a_m;c)$ does not exist. 
\end{cor}
 

\begin{thm} \label{reduction}
Let $a_1,\ldots,a_m$ be nonzero integers, $c \in \Z$, $r>1$ and $\l \ge 1$. If $\sum_{i=1}^m a_i=\s$, then
\[ {\tt Rad}_r (a_1,\ldots,a_m;\l c+\s) \le 1 + \l \big( {\tt Rad}_r (a_1,\ldots,a_m;c+\s) - 1 \big).  \]
\end{thm}

\begin{Pf}
Let $N={\tt Rad}_r (a_1,\ldots,a_m;c+\s)$. Then any $r$-colouring of the integers in $[1,N]$ must contain a monochromatic solution to $\sum_{i=1}^m a_ix_i=c+\s$. We extend any such $r$-colouring $\chi$ to $[0,\l(N-1)]$ by first assigning the colour of $x$ in $[1,N]$ to $x-1$. Thus we have an $r$-colouring of the integers in $[0,N-1]$. Then we assign $x,y$ in $[0,\l(N-1)]$ the same colour if and only if $x \equiv y\mod{N-1}$. Since solutions $(x_1,\ldots,x_m)$ to $\sum_{i=1}^m a_ix_i=c+\s$ give rise to solutions $\big(\l(x_1-1),\ldots,\l(x_m-1)\big)$ to $\sum_{i=1}^m a_i\big(\l (x_i-1)+1\big)=\l c+\s$, the result follows from the argument at the beginning of this section. 
\end{Pf}

In this paper, we study the Rado numbers for the equation 
\begin{equation} \label{spleqn}
\sum_{i=1}^{m-1}a_ix_i - x_m = c
\end{equation}
in the special case where $\{a_1,\ldots,a_{m-1}\}$ is $2$-distributable or $3$-distributable and $r=2$, for any $c \in \Z$. By Corollary \ref{dne}, ${\tt Rad}_2 \big(a_1,\ldots,a_{m-1},-1;c\big)$ exists if and only if either $c$ is even or $\S$ is even. Throughout the rest of this paper, we assume that at least one of $c$, $\S$ is even. 

By assigning the colour of $x_i$ in the solution of Eq.(\ref{spleqn}) to $x_i-1$, we note that this is equivalent to determining the smallest positive integer $R$ for which every $r$-colouring of $[0,R-1]$ contains a monochromatic solution to 
\begin{equation} \label{revspleqn}
\sum_{i=1}^{m-1} a_ix_i - x_m = c-(\S-1),
\end{equation}
where $\S=\sum_{i=1}^{m-1} a_i$. 

For brevity, we write 
\[ {\tt Rad}_r \big(c\big) \]
for the $r$-colour Rado number ${\tt Rad}_r \big(a_1,\ldots,a_{m-1},-1;c\big)$, and 
\[ {\tt Rad}_{r;t} \big(c\big) \]
for the $r$-colour Rado number ${\tt Rad}_r \big(a_1,\ldots,a_{m-1},-1;c\big)$ for $t$-distributable sets $\{a_1,\ldots,a_{m-1}\}$. 

\begin{thm} \label{c<S-1}
Let $a_1,\ldots,a_{m-1}$ be a set of positive integers, with $\sum_{i=1}^{m-1} a_i=\S$.  
\begin{itemize}
\item[{\rm (i)}] 
For any $r>1$, 
\[ {\tt Rad}_r \big(\S-1\big) = 1. \]
\item[{\rm (ii)}] 
If $c\,\S$ is even and $c<\S-1$, then 
\[ {\tt Rad}_2 \big(c\big) \ge (\S-1-c)(\S+2) + 1. \]
\item[{\rm (iii)}] 
If $c\,\S$ is even, $c<\S-1$, and if the set $\{a_1,\ldots,a_{m-1}\}$ is $2$-distributable, then 
\[ {\tt Rad}_{2;2} \big(c\big) = (\S-1-c)(\S+2) + 1. \]
\end{itemize}
\end{thm}

\begin{Pf}
\begin{itemize}
\item[{\rm (i)}]
If $c=\S-1$, then Eq.(\ref{spleqn}) admits the solution $x_1=\cdots=x_m=1$. Hence every $r$-colouring of the set $\{1\}$ admits a monochromatic solution of Eq.(\ref{spleqn}). Thus ${\tt Rad}_r \big(\S-1\big)=1$ for each $r>1$.   
\item[{\rm (ii)}]
By Corollary \ref{dne}, ${\tt Rad}_2 \big(c\big)$ exists if and only if $c\,\S$ is even. For the rest of the proof, we also suppose $c<\S-1$. Let $\D: [0,(\S-1-c)(\S+2)-1] \imp \{0,1\}$ be defined by 
\[ {\D}(x) = \begin{cases}
                       0 & \mbox{ if $x \in [0,\S-2-c] \cup [(\S-1-c)(\S+1), (\S-1-c)(\S+2)-1]$}; \\
                       1 & \mbox{ if $x \in [\S-1-c,(\S-1-c)(\S+1)-1]$}.
                      \end{cases}
\]
We claim that $\D$ provides a valid $2$-colouring of $[0,(\S-1-c)(\S+2)-1]$ with respect to Eq.(\ref{revspleqn}). 

Suppose ${\D}(x_i)=1$ for $i \in \{1,\ldots,m-1\}$. Then
\[ x_m = \sum_{i=1}^{m-1} a_ix_i + \S-1-c  \ge \S(\S-1-c) + (\S-1-c) = (\S+1)(\S-1-c). \] 
Hence ${\D}(x_m)=0$. Therefore, we must have ${\D}(x_i)=0$ for $i \in \{1,\ldots,m\}$. 

If $x_i \in [0,\S-2-c]$ for $i \in \{1,\ldots,m-1\}$, then
\[ x_m = \sum_{i=1}^{m-1} a_ix_i + \S-1-c  \ge \S-1-c, \]
and 
\[ x_m = \sum_{i=1}^{m-1} a_ix_i + \S-1-c  \le \S(\S-2-c) + (\S-1-c)  < (\S-1-c)(\S+1). \] 
Therefore, $x_i \in [(\S-1-c)(\S+1), (\S-1-c)(\S+2)-1]$ for at least one $i \in \{1,\ldots,m-1\}$. Now 
\[ x_m = \sum_{i=1}^{m-1} a_ix_i + \S-1-c  \ge (\S-1-c)(\S+1) + (\S-1-c) = (\S-1-c)(\S+2), \]
so that $x_m$ is outside the domain of $\D$. Hence ${\tt Rad}_2 \big(c\big) \ge (\S-1-c)(\S+2)+1$. 
\item[{\rm (iii)}]
Suppose $\{a_1,\ldots,a_{m-1}\}$ is $2$-distributable. For any ${\chi}: [0,(\S-1-c)(\S+2)] \imp \{0,1\}$, we show that there exists a monochromatic solution to Eq.(\ref{revspleqn}) under $\chi$. Together with the result in part (ii), this would complete the proof of part (iii). We first resolve the case where $c=\S-2$; under our assumption, both $c$ and $\S$ are even. 

Let $c=\S-2$, and let ${\chi}: [0,\S+2] \imp \{0,1\}$. Without loss of generality, let ${\chi}(0)=0$. If ${\chi}(1)=0$, then $x_1=\cdots=x_{m-1}=0$, $x_m=1$ provides a monochromatic solution to Eq.(\ref{revspleqn}). Hence ${\chi}(1)=1$. 

If ${\chi}(\S+1)=1$, then $x_1=\cdots=x_{m-1}=1$, $x_m=\S+1$ provides a monochromatic solution to Eq.(\ref{revspleqn}). Hence ${\chi}(\S+1)=0$. 

Since $\{a_1,\ldots,a_{m-1}\}$ is $2$-distributable with $\sum_{i=1}^{m-1} a_i=\S$, distribute the sets according to ${\s}_1=1$ and ${\s}_2=\S-1$. We may assume that $a_1=1$, after relabelling if necessary. If ${\chi}(\S+2)=0$, then $x_1=\S+1$, $x_2=\cdots=x_{m-1}=0$ provides a monochromatic solution to Eq.(\ref{revspleqn}). Hence ${\chi}(\S+2)=1$. 

We now distribute the set $\{a_1,\ldots,a_{m-1}\}$ according to ${\s}_1={\s}_2=\S/2$. If ${\chi}(2)=0$, then $x_i=0$ for the $a_i$'s in the first collection and $x_j=2$ for the $a_j$'s in the second collection provides a monochromatic solution to Eq.(\ref{revspleqn}) since ${\chi}(\S+1)=0$. Hence ${\chi}(2)=1$. 

Since $\{a_1,\ldots,a_{m-1}\}$ is $2$-distributable, we may assume that at least one $a_i=1$; after relabelling we may assume $a_1=1$. Now $x_1=2$, $x_2=\cdots=x_m=1$ provides a monochromatic solution to Eq.(\ref{revspleqn}) since ${\chi}(1)={\chi}(2)={\chi}(\S+2)=1$. 

Therefore any $2$-colouring of the integers in $[0,\S+2]$ admits a monochromatic solution to $\sum_{i=1}^{m-1} a_ix_i + 1=x_m$, which is Eq.(\ref{revspleqn}) for $c=\S-2$. Hence ${\tt Rad}_{2;2} \big(\S-2\big) \le \S+3$. 
  
Now let $c<\S-1$. With $c=-1$ and $\l=(\S-1)-c$, Theorem \ref{reduction} gives
\begin{eqnarray*}
{\tt Rad}_{2;2} \big(\s-(\S-1-c)\big) \le 1 + (\S-1-c) \big( {\tt Rad}_{2;2} (\s-1) - 1 \big),
\end{eqnarray*}
or that 
\begin{eqnarray*}
{\tt Rad}_{2;2} \big(c\big) \le 1 + (\S-1-c)(\S+2). 
\end{eqnarray*}
\end{itemize}
\end{Pf}


\begin{thm} \label{[S,2S-1]}
Let $a_1,\ldots,a_{m-1}$ be a set of positive integers, with $\sum_{i=1}^{m-1} a_i=\S$. Let $\S \le c \le 2\S-1$ and let $c\,\S$ be even.  
\begin{itemize}
\item[{\rm (i)}] 
\[ {\tt Rad}_r \big(2\S-2\big) = 2. \]
\item[{\rm (ii)}] 
Suppose $\{a_1,\ldots,a_{m-1}\}$ is $2$-distributable. If $c \ne 2\S-2$, then 
\[ {\tt Rad}_{2;2} \big(c\big) \le \S+1 \;\;\text{and}\;\; {\tt Rad}_{2;2} \big(2\S-1\big) = 3. \]
\end{itemize}
\end{thm}

\begin{Pf}
\begin{itemize}
\item[{\rm (i)}] 
If $c=2\S-2$, then Eq.(\ref{spleqn}) admits the solution $x_1=\cdots=x_m=2$. Hence every $r$-colouring of the set $\{1,2\}$ admits a monochromatic solution of Eq.(\ref{spleqn}). Since $x_1=\cdots=x_m=1$ is not a solution to Eq.(\ref{spleqn}),  it follows that ${\tt Rad}_r \big(2\S-2\big)=2$ for each $r>1$.  

\item[{\rm (ii)}] 
Suppose $\{a_1,\ldots,a_{m-1}\}$ is $2$-distributable, and let $r=2$. 

Suppose $c=2\S-1$; since $c$ is odd, $\S$ must be even. For any ${\chi}: [1,3] \imp \{0,1\}$, we show that there exists a monochromatic solution to Eq.(\ref{spleqn}) under $\chi$. 

If ${\chi}(1)={\chi}(2)$, then Eq.(\ref{spleqn}) admits the solution $x_1=\cdots=x_{m-1}=2$, $x_m=1$. 

Henceforth suppose ${\chi}(1) \ne {\chi}(2)$. If ${\chi}(3)={\chi}(1)$, distribute the set $\{a_1,\ldots,a_{m-1}\}$ according to ${\s}_1={\s}_2=\S/2$. 
Then $x_i=1$ for $a_i$'s in the first collection, $x_i=3$ for $a_i$'s in the second collection, and $x_m=1$ provides a monochromatic solution to Eq.(\ref{spleqn}). If ${\chi}(3)={\chi}(2)$, distribute the set $\{a_1,\ldots,a_{m-1}\}$ according to ${\s}_1=\S-1$ and ${\s}_2=1$. Then $x_i=2$ for $a_i$'s in the first collection, $x_i=3$ for $a_i$ $(=1)$ in the second collection, and $x_m=2$ provides a monochromatic solution to Eq.(\ref{spleqn}). 

This proves our claim, and shows that ${\tt Rad}_{2;2} \big(c\big) \le 3$ in this case. 

To prove that ${\tt Rad}_{2;2} \big(c\big)>2$ in this case, observe that the $2$-colouring of $[1,2]$ with ${\chi}(1) \ne {\chi}(2)$ is a valid colouring. This completes the proof in the special case $c=2\S-1$. 

It remains to show that ${\tt Rad}_{2;2} \big(c\big) \le \S+1$ for $c \in [\S,2\S-3]$. Write $c=2(\S-1)-k$, with $k \in [1,\S-2]$. By assigning the colour of $x_i$ in the solution of Eq.(\ref{spleqn}) to $x_i-1$, we show that any $2$-colouring of $[0,\S]$ admits a monochromatic solution to 
\begin{equation} \label{revspleqn_1}
\sum_{i=1}^{m-1} a_ix_i - x_m = \S-1-k,
\end{equation}
for $k \in [1,\S-2]$. 

Let ${\chi}: [0,\S] \imp \{0,1\}$, and assume ${\chi}(1)=0$ without loss of generality. 

If ${\chi}(k+1)=0$, then $x_1=\cdots=x_{m-1}=1$, $x_m=k+1$ provides a monochromatic solution to Eq.(\ref{revspleqn_1}). Therefore we may assume ${\chi}(k+1)=1$. 

Suppose ${\chi}(0)=0$. Distribute the set $\{a_1,\ldots,a_{m-1}\}$ according to ${\s}_1=k+1$ and ${\s}_2=\S-(k+1)$. Then $x_i=0$ for the $a_i$'s in the first collection, $x_i=1$ for the $a_i$'s in the second collection, and $x_m=1$ provides a monochromatic solution to Eq.(\ref{revspleqn_1}). Therefore we may assume ${\chi}(0)=1$. 

Suppose ${\chi}(\S-1-k)=1$. Distribute the set $\{a_1,\ldots,a_{m-1}\}$ according to ${\s}_1=\S-1$ and ${\s}_2=1$. Then $x_i=0$ for the $a_i$'s in the first collection, $x_i=\S-1-k$ for the $a_i$ $(=1)$ in the second collection, and $x_m=0$ provides a monochromatic solution to Eq.(\ref{revspleqn_1}). Therefore we may assume ${\chi}(\S-1-k)=0$. 
 
Suppose ${\chi}(\S)=1$. Distribute the set $\{a_1,\ldots,a_{m-1}\}$ according to ${\s}_1=\S-1$ and ${\s}_2=1$. Then $x_i=0$ for the $a_i$'s in the first collection, $x_i=\S$ for the $a_i$ $(=1)$ in the second collection, and $x_m=k+1$ provides a monochromatic solution to Eq.(\ref{revspleqn_1}). Therefore we may assume ${\chi}(\S)=0$. 

Suppose ${\chi}(\S-1)=0$. Distribute the set $\{a_1,\ldots,a_{m-1}\}$ according to ${\s}_1=\S-1$ and ${\s}_2=1$. Then $x_i=1$ for the $a_i$'s in the first collection, $x_i=\S-1-k$ for the $a_i$ $(=1)$ in the second collection, and $x_m=\S-1$ provides a monochromatic solution to Eq.(\ref{revspleqn_1}). Therefore we may assume ${\chi}(\S-1)=1$. 

Suppose $c$ is even; then $k$ is even. Distribute the set $\{a_1,\ldots,a_{m-1}\}$ according to ${\s}_1=\frac{k}{2}+1$ and ${\s}_2=\S-\big(\frac{k}{2}+1\big)$. Then $x_i=0$ for the $a_i$'s in the first collection, $x_i=2$ for the $a_i$'s in the second collection, and $x_m=\S-1$ provides a monochromatic solution to Eq.(\ref{revspleqn_1}). Therefore we may assume ${\chi}(2)=0$ in this case. 

Now suppose $c$ is odd; then $k$ is odd and $\S$ is even. Distribute the set $\{a_1,\ldots,a_{m-1}\}$ according to ${\s}_1=(\S+1+k)/2$ and ${\s}_2=(\S-1-k)/2$. Then $x_i=0$ for the $a_i$'s in the first collection, $x_i=2$ for the $a_i$'s in the second collection, and $x_m=0$ provides a monochromatic solution to Eq.(\ref{revspleqn_1}). Therefore we may also assume ${\chi}(2)=0$ in this case. 

In both cases, we must have ${\chi}(2)=0$. Distribute the set $\{a_1,\ldots,a_{m-1}\}$ according to ${\s}_1=k+1$ and ${\s}_2=\S-(k+1)$. Then $x_i=1$ for the $a_i$'s in the first collection, $x_i=2$ for the $a_i$'s in the second collection, and $x_m=\S$ provides a monochromatic solution to Eq.(\ref{revspleqn_1}). 
\end{itemize}
\end{Pf}


\begin{thm} \label{lambda_mu_1}
Let $a_1,\ldots,a_{m-1}$ be a set of positive integers, with $\sum_{i=1}^{m-1} a_i=\S$. Let $c \in \bigcup_{3 \le \l \le \S} \big[(\l-1)\S,\l\S-\l\big)$, and let $c\,\S$ be even. Then $c=\l(\S-1)-\mu$, with $3 \le \l \le \S$, $0 \le \mu \le \S-\l$, and  
\begin{itemize}
\item[{\rm (i)}] 
\[ {\tt Rad}_2 \big(c\big) \ge \l+\mu. \]
\item[{\rm (ii)}] 
For any $r>1$, 
\[ {\tt Rad}_r \big(\l(\S-1)\big) = \l. \]
\item[{\rm (iii)}] 
If the set $\{a_1,\ldots,a_{m-1}\}$ is $2$-distributable, then 
\[ {\tt Rad}_{2;2} \big(\l(\S-1)-1\big) = \l+1. \]
\item[{\rm (iv)}] 
If the set $\{a_1,\ldots,a_{m-1}\}$ is $3$-distributable and $\mu \ne 0,1$, then 
\[ {\tt Rad}_{2;3} \big(c\big) = \l+\mu. \]
\end{itemize}
\end{thm}

\begin{Pf}
By Corollary \ref{dne}, ${\tt Rad}_2 \big(c\big)$ exists if and only if $c\,\S$ is even. Suppose $c=\l(\S-1)-\mu$, with $3 \le \l \le \S$, $0 \le \mu \le \S-\l$. 
\begin{itemize}
\item[{\rm (i)}] 
Let $\D: [0,\l+\mu-2] \imp \{0,1\}$ be defined by 
\[ {\D}(x) = \begin{cases}
                       0 & \mbox{ if $x \in [0,\l-2]$}; \\
                       1 & \mbox{ if $x \in [\l-1,\l+\mu-2]$}.
                      \end{cases}
\]
We claim that $\D$ provides a valid $2$-colouring of $[0,\l+\mu-2]$ with respect to Eq.(\ref{revspleqn}). 

Suppose ${\D}(x_i)=0$ for $i \in \{1,\ldots,m-1\}$. 
\[ x_m = \sum_{i=1}^{m-1} a_ix_i - (\l-1)(\S-1) + \mu \le \S(\l-2) - (\l-1)(\S-1) + \mu \le -1, \] 
so that $x_m$ is outside the domain of $\D$. Therefore, we must have ${\D}(x_i)=1$ for $i \in \{1,\ldots,m\}$. But then 
\[ x_m = \sum_{i=1}^{m-1} a_ix_i - (\l-1)(\S-1) + \mu \ge \S(\l-1) - (\l-1)(\S-1) + \mu \ge \l+\mu-1. \] 
Hence ${\tt Rad}_2 \big(c\big) \ge \l+\mu$. 

\item[{\rm (ii)}] 
Let $r>1$. By part (i), 
\[ {\tt Rad}_r \big(\l(\S-1)\big) \ge {\tt Rad}_2 \big(\l(\S-1)\big) \ge \l. \] 

Since Eq.(\ref{spleqn}) admits the solution $x_1=\cdots=x_m=\l$, every $r$-colouring of $[1,\l]$ admits a monochromatic solution of Eq.(\ref{spleqn}). Therefore ${\tt Rad}_r \big(\l(\S-1)\big) \le \l$ for each $r>1$.  

\item[{\rm (iii)}] 
Let $c=\l(\S-1)-\mu$ with $3 \le \l \le \S$, $1 \le \mu \le \S-\l$. By part (i), it suffices to show that 
\[ {\tt Rad}_{2;2} \big(c\big) \le \l+\mu \] 
when the set $\{a_1,\ldots,a_{m-1}\}$ is $2$-distributable for $\mu=1$, and $3$-distributable for $\mu>1$. 

By assigning the colour of $x_i$ in the solution of Eq.(\ref{spleqn}) to $x_i-1$, we show that any $2$-colouring of $[0,\l+\mu-1]$ admits a monochromatic solution to 
\begin{equation} \label{revspleqn_2}
\sum_{i=1}^{m-1} a_ix_i - x_m = (\l-1)(\S-1)-\mu.
\end{equation}

Let ${\chi}: [0,\l+\mu-1] \imp \{0,1\}$, and assume ${\chi}(\l-1)=0$ without loss of generality. Assume that the set $\{a_1,\ldots,a_{m-1}\}$ is $2$-distributable.  

If ${\chi}(\l-1+\mu)=0$, then $x_1=\cdots=x_{m-1}=\l-1$, $x_m=\l-1+\mu$ provides a monochromatic solution to Eq.(\ref{revspleqn_2}). Therefore we may assume ${\chi}(\l-1+\mu)=1$. 

Suppose ${\chi}(\l-2)=0$. Distribute the set $\{a_1,\ldots,a_{m-1}\}$ according to ${\s}_1=\mu$ and ${\s}_2=\S-\mu$. Then $x_i=\l-2$ for the $a_i$'s in the first collection, $x_i=\l-1$ for the $a_i$'s in the second collection, and $x_m=\l-1$ provides a monochromatic solution to Eq.(\ref{revspleqn_2}). Therefore we may assume ${\chi}(\l-2)=1$. 

Suppose $c$ is odd; then $\mu$ is odd and $\S$ is even. Distribute the set $\{a_1,\ldots,a_{m-1}\}$ according to ${\s}_1={\s}_2=\S/2$. Then $x_i=\l$ for the $a_i$'s in the first collection, $x_i=\l-2$ for the $a_i$'s in the second collection, and $x_m=\l-1+\mu$ provides a monochromatic solution to Eq.(\ref{revspleqn_2}). Therefore we may also assume ${\chi}(\l)=0$ in this case. 

Then $x_1=\cdots=x_{m-1}=\l-1$ and $x_m=\l-1+\mu$ provides a solution to Eq.(\ref{revspleqn_2}). However, ${\chi}(\l-1)=0$, whereas ${\chi}(\l-1+\mu)=1$ for $\mu>1$ and ${\chi}(\l)=0$. Hence we have a monochromatic solution to Eq.(\ref{revspleqn_2}) for $\mu=1$.   
 
\item[{\rm (iv)}] 
For the rest of this proof, assume that $\mu>1$. If $c$ is odd, again distribute the set $\{a_1,\ldots,a_{m-1}\}$ according to ${\s}_1=(\mu-1)/2$ and ${\s}_2=\S-(\mu-1)/2$. Then $x_i=\l-3$ for the $a_i$'s in the first collection, $x_i=\l-1$ for the $a_i$'s in the second collection, and $x_m=\l$ provides a monochromatic solution to Eq.(\ref{revspleqn_2}). Therefore we may also assume ${\chi}(\l-3)=1$ in this case. 

Next suppose $c$ is even; then $\mu$ is even. Distribute the set $\{a_1,\ldots,a_{m-1}\}$ according to ${\s}_1=\mu/2$ and ${\s}_2=\S-\mu/2$. Then $x_i=\l-3$ for the $a_i$'s in the first collection, $x_i=\l-1$ for the $a_i$'s in the second collection, and $x_m=\l-1$ provides a monochromatic solution to Eq.(\ref{revspleqn_2}). Therefore we may assume ${\chi}(\l-3)=1$ in this case as well. 

Henceforth assume that the set $\{a_1,\ldots,a_{m-1}\}$ is $3$-distributable. Distribute the set $\{a_1,\ldots,a_{m-1}\}$ according to ${\s}_1=\big\lc\frac{\S}{\mu+1}\big\rc$, ${\s}_2=\big\lc\frac{\S}{\mu+1}\big\rc (\mu+1)-\S$, and ${\s}_3=2\S-\big\lc\frac{\S}{\mu+1}\big\rc (\mu+2)$. Then $x_i=\l-1+\mu$ for the $a_i$'s in the first collection, $x_i=\l-3$ for the $a_i$'s in the second collection, $x_i=\l-2$ for the $a_i$'s in the third collection, together with 
\begin{eqnarray*}
x_m & = & \sum_{i=1}^{m-1} a_ix_i - (\l-1)(\S-1) + \mu \\
& = & (\l-1+\mu) \left\lc\ts\frac{\S}{\mu+1}\right\rc + (\l-3) \left( \left\lc\ts\frac{\S}{\mu+1}\right\rc (\mu+1) - \S \right) + (\l-2) \left( 2\S - 
\left\lc\ts\frac{\S}{\mu+1}\right\rc (\mu+2) \right) \\  
& & - (\l-1)(\S-1) + \mu \\
& = & \l - 1 + \mu
\end{eqnarray*} 
provides a monochromatic solution to Eq.(\ref{revspleqn_2}). 
\end{itemize}
\end{Pf}


\begin{thm} \label{lambda_mu_2}
Let $a_1,\ldots,a_{m-1}$ be a $3$-distributable set of positive integers, with $\sum_{i=1}^{m-1} a_i=\S$. Let $c \in \bigcup_{3 \le \l \le \big\lc\frac{\S+1}{2}\big\rc} \big(\l(\S-1),\l\S\big)$, and let $c\,\S$ be even. If $c=\l\S-\mu$, with $3 \le \l \le \big\lc\frac{\S+1}{2}\big\rc$, $1 \le \mu \le \l-1$, then 
\[ {\tt Rad}_{2;3} \big(c\big) \le 2\l-\mu. \]
\end{thm}

\begin{Pf}
Let $c=\l\S-\mu$, with $3 \le \l \le \big\lc\frac{\S+1}{2}\big\rc$ and $1 \le \mu \le \l-1$. By assigning the colour of $x_i$ in the solution of Eq.(\ref{spleqn}) to $x_i-1$, we show that any $2$-colouring of $[0,2\l-\mu-1]$ admits a monochromatic solution to 
\begin{equation} \label{revspleqn_3}
\sum_{i=1}^{m-1} a_ix_i - x_m = (\l-1)\S - (\mu-1).
\end{equation}

Let ${\chi}: [0,2\l-\mu-1] \imp \{0,1\}$, and assume ${\chi}(\l-1)=0$ without loss of generality. Assume that the set $\{a_1,\ldots,a_{m-1}\}$ is $3$-distributable.  

If ${\chi}(\mu-1)=0$, then $x_1=\cdots=x_{m-1}=\l-1$, $x_m=\mu-1$ provides a monochromatic solution to Eq.(\ref{revspleqn_3}). Therefore we may assume ${\chi}(\mu-1)=1$. 

Suppose ${\chi}(2\l-\mu-1)=0$. Distribute the set $\{a_1,\ldots,a_{m-1}\}$ according to ${\s}_1=\S-1$ and ${\s}_2=1$. Then $x_i=\l-1$ for the $a_i$'s in the first collection, $x_i=2\l-\mu-1$ for the $a_i$ (=1) in the second collection, and $x_m=\l-1$, provides a monochromatic solution to Eq.(\ref{revspleqn_3}). Therefore we may assume ${\chi}(2\l-\mu-1)=1$. 

Suppose $\S$ is even. Distribute the set $\{a_1,\ldots,a_{m-1}\}$ according to ${\s}_1=\S/2$ and ${\s}_2=\S/2$. Then $x_i=2\l-\mu-1$ for the $a_i$'s in the first collection, $x_i=\mu-1$ for the $a_i$'s in the second collection, and $x_m=\mu-1$, provides a monochromatic solution to Eq.(\ref{revspleqn_3}). 

Now suppose $\S$ is odd, so that $c=\l\S-\mu$ is even. Hence $\l-\mu$ is even. Distribute the set $\{a_1,\ldots,a_{m-1}\}$ according to ${\s}_1=\l-\mu$ and ${\s}_2=\S-(\l-\mu)$. Then $x_i=\l$ for the $a_i$'s in the first collection, $x_i=\l-1$ for the $a_i$'s in the second collection, and $x_m=\l-1$, provides a solution to Eq.(\ref{revspleqn_3}). Therefore we may assume ${\chi}(\l)=1$.

Distribute the set $\{a_1,\ldots,a_{m-1}\}$ according to ${\s}_1=\l-\mu$, ${\s}_2=(\S-\l+\mu+1)/2$ and ${\s}_3=(\S-\l+\mu-1)/2$. Then $x_i=\l$ for the $a_i$'s in the first collection, $x_i=\mu-1$ for the $a_i$'s in the second collection, $x_i=2\l-\mu-1$ for the $a_i$'s in the third collection, and 
$x_m=\mu-1$, provides a solution to Eq.(\ref{revspleqn_3}). 
\end{Pf}


\begin{thm} \label{>S(S-1)}
Let $\{a_1,\ldots,a_{m-1}\}$ be a set of positive integers, with $\sum_{i=1}^{m-1} a_i=\S$. Let $c>\S-1$, and let $c\,\S$ be even.  
\begin{itemize}
\item[{\rm (i)}] 
\[ {\tt Rad}_2 \big(c\big) \ge \left\lc\frac{1+c(\S+2)}{{\S}^2+\S-1}\right\rc. \]
\item[{\rm (ii)}] 
If the set $\{a_1,\ldots,a_{m-1}\}$ is $2$-distributable and $c>\S(\S-1)$, then 
\[ {\tt Rad}_{2;2} \big(c\big) = \left\lc\frac{1+c(\S+2)}{{\S}^2+\S-1}\right\rc. \]
\end{itemize}
\end{thm}

\begin{Pf}
For convenience, we set 
\[ \T = \left\lc\frac{1+c(\S+2)}{{\S}^2+\S-1}\right\rc - 1, \]
and show that 
\begin{equation} \label{L_ineq}
{\S}^2 \T - c(\S+1) < \S \T - c < \T. 
\end{equation}
Both inequalities are both equivalent to 
\[ \T < \frac{c}{\S-1}. \]
From the definition of $\T$, since $c>\S-1$, we have 
\begin{equation} \label{T_upper}
\T \le \frac{1+c(\S+2)}{{\S}^2+\S-1} = \frac{1+c(\S+2)}{1+(\S-1)(\S+2)} < \frac{c}{\S-1}. 
\end{equation}
Thus both inequalities in Eq.(\ref{L_ineq}) hold. 

We have
\begin{eqnarray*} 
\S \T - c & \ge & \S \left( \frac{1+c(\S+2)}{{\S}^2+\S-1} - 1 \right) - c \\
& = & \S \left(\frac{1+c(\S+2)}{1+(\S-1)(\S+2)} - 1 \right) - c \\
& = &  \frac{\big(c-(\S-1)\big)\S(\S+2)}{1+(\S-1)(\S+2)} - c \\
& = & \frac{c(\S+1)-\S(\S-1)(\S+2)}{1+(\S-1)(\S+2)} \\
& > & 0 \;\;\text{if}\;\; c>(\S+1)(\S-1).
\end{eqnarray*}

\begin{itemize}
\item[{\rm (i)}]
Let $\D: [1,\T] \imp \{0,1\}$ be defined by 
\[ {\D}(x) = \begin{cases}
                       0 & \mbox{ if $x \in \big( \max \{0,{\S}^2 \T - c(\S+1)\},\S \T - c \,\big]$}; \\
                       1 & \mbox{ otherwise}.
                      \end{cases}
\]
We claim that $\D$ provides a valid $2$-colouring of $[1,\T]$ with respect to Eq.(\ref{spleqn}). 

Suppose ${\D}(x_i)=0$ for $i \in \{1,\ldots,m-1\}$. Then
\[ x_m = \sum_{i=1}^{m-1} a_i x_i - c  \le \S(\S \T - c) - c = {\S}^2 \T - c(\S+1). \] 
Hence ${\D}(x_m)=1$. Therefore, we must have ${\D}(x_i)=1$ for $i \in \{1,\ldots,m\}$. 

If $x_i>\S \T - c$ for $i \in \{1,\ldots,m-1\}$, then
\[ x_m = \sum_{i=1}^{m-1} a_i x_i - c  > \S(\S \T - c) - c = {\S}^2 \T - c(\S+1), \]
and 
\[ x_m = \sum_{i=1}^{m-1} a_i x_i - c  \le \S \T - c. \] 
Therefore, $x_i \in \big[1,{\S}^2 \T - c(\S+1)\big]$ for at least one $i \in \{1,\ldots,m-1\}$. Now 
\begin{eqnarray*} 
x_m & = & \sum_{i=1}^{m-1} a_i x_i - c  \\
& \le & {\S}^2 \T - c(\S+1) + (\S-1)\T - c \\
& = & ({\S}^2+\S-1) \T - c(\S+2) \\
& < & ({\S}^2+\S-1) \cdot \frac{1+c(\S+2)}{{\S}^2+\S-1} - c(\S+2) \\
& = & 1, 
\end{eqnarray*}
so that $x_m$ is outside the domain of $\D$. 

We note that ${\S}^2 \T - c(\S+1) < \S \T - c$ for $c>\S-1$, and further that $\S \T-c>0$ if $c>(\S+1)(\S-1)$. Thus $\D$ provides a valid $2$-colouring for $c>(\S+1)(\S-1)$. For $c \in [\S,{\S}^2)$, it may be the case that $\S \T-c<1$, in which case all integers in the interval $[1,\T]$ are coloured $1$. Since $x_m=\sum_{i=1}^{m-1} a_i x_i - c \le \S \T-c$, $\D$ provides a valid $1$-colouring if $\S \T-c<1$. Therefore ${\tt Rad}_1 \big(c\big)>\T$ in such cases, and ${\tt Rad}_2 \big(c\big)>\T$ in any case. 

\item[{\rm (ii)}]
Suppose that the set $\{a_1,\ldots,a_{m-1}\}$ is $2$-distributable, and $c>\S(\S-1)$. By part (i), it suffices to prove that 
\[ {\tt Rad}_{2;2} \big(c\big) \le \left\lc\frac{1+c(\S+2)}{{\S}^2+\S-1}\right\rc = \T+1. \] 

We have
\[ \T \ge \frac{\big(c-(\S-1)\big)(\S+2)}{{\S}^2+\S-1} \ge \frac{\big(1+\S(\S-1)-(\S-1)\big)(\S+2)}{{\S}^2+\S-1} = \S-1 + \frac{3}{{\S}^2+\S-1}. \]
Hence $T \ge \S$ when $c>\S(\S-1)$. 

Let $\chi: [1,\T+1] \imp \{0,1\}$ be any $2$-colouring of the integers in the interval $[1,\T+1]$. Consider the complimentary colouring $\ov{\chi}: [1,\T+1] \imp \{0,1\}$ given by 
\[ \ov{\chi}(x) = \chi(\T+2-x). \]
Then monochromatic solutions to Eq.(\ref{spleqn}) under $\chi$ correspond to monochromatic solutions to 
\begin{equation} \label{compspleqn}
\sum_{i=1}^{m-1} a_ix_i - x_m = (\S-1)(\T+2) - c 
\end{equation}
under $\ov{\chi}$. 

From Eq.(\ref{T_upper}), 
\[  (\S-1)(\T+2) - c < (\S-1)\left(\frac{c}{\S-1}+2\right) - c = 2(\S-1). \]

Thus we have $c^{\pr}=(\S-1)(\T+2)-c<2(\S-1)$ for $c>\S(\S-1)$. If $c^{\pr}<\S-1$, every $2$-colouring of $[1,(\S-1-c^{\pr})(\S+2)+1]$ admits a monochromatic solution to Eq.(\ref{compspleqn}) by Theorem \ref{c<S-1}. Now 
\begin{eqnarray*}
(\S-1-c^{\pr})(\S+2)+1 & = & \big(c - (\S-1)(\T+1)\big)(\S+2)+1 \\ 
& = & 1 + c(\S+2) - \big(1+(\S-1)(\S+2)\big)(\T+1) + \T+1 \\
& \le & \T+1. 
\end{eqnarray*}
Hence every $2$-colouring of $[1,\T+1]$ also admits a monochromatic solution to Eq.(\ref{compspleqn}) in this case. If $c^{\pr} \in \big[\S-1,2(\S-1)\big)$, then every $2$-colouring of $[1,\S+1]$ admits a monochromatic solution to Eq.(\ref{compspleqn}) by Theorem \ref{[S,2S-1]}. Since $\S \le \T$, every $2$-colouring of $[1,\T+1]$ also admits a monochromatic solution to Eq.(\ref{compspleqn}). 
\end{itemize}
\end{Pf}


\begin{thm} \label{S+1}
Let $\{a_1,\ldots,a_{m-1}\}$ be a $3$-distributable set of positive integers, with $\sum_{i=1}^{m-1} a_i=\S$. If $c \in \big(\S-1,\S(\S-1)\big)$ and $c\,\S$ is  even, then  
\[ {\tt Rad}_{2;3} \big(c\big) \le \S+1. \]
\end{thm}

\begin{Pf}
The result holds for $2$-distributable sets $\{a_1,\ldots,a_{m-1}\}$ when $c \in (\S-1,2\S-1]$ by Theorem \ref{[S,2S-1]}. 

The range of $c$ in Theorems \ref{lambda_mu_1} and \ref{lambda_mu_2} together cover $\bigcup_{3 \le \l \le \lc\frac{\S+1}{2}\rc}  
[(\l-1)\S,\l\S)=[2\S,\lc\frac{\S+1}{2}\rc \S)$. In the cases covered by Theorem \ref{lambda_mu_1}, the Rado number equals $\l+\mu$, which is at most $\S$. In the cases covered by Theorem \ref{lambda_mu_2}, the Rado number is at most $2\l-\mu$, and this is at most $2 \frac{\S+2}{2}-1=\S+1$.  Therefore for $3$-distributable sets of positive integers $a_1,\ldots,a_{m-1}$ and for $c \in [2\S,\lc\frac{\S+1}{2}\rc \S)$, we have  
\[ {\tt Rad}_{2;3} \big(c\big) \le \S+1. \]   

Let $\chi: [1,\S+1] \imp \{0,1\}$ be any $2$-colouring of the integers in the interval $[1,\S+1]$. Consider the complimentary colouring $\ov{\chi}: [1,\S+1] \imp \{0,1\}$ given by 
\[ \ov{\chi}(x) = \chi(\S+2-x). \]
Then monochromatic solutions to Eq.(\ref{spleqn}) under $\chi$ correspond to monochromatic solutions to 
\begin{equation} \label{compspleqn}
\sum_{i=1}^{m-1} a_ix_i - x_m = (\S-1)(\S+2) - c 
\end{equation}
under $\ov{\chi}$. 

Note that 
\[ \left\lc\ts\frac{\S+1}{2}\right\rc \S \le c < \S(\S-1) \;\;{\rm implies}\;\; 2(\S-1) < (\S-1)(\S+2)-c \le (\S-1)(\S+2)-\left\lc\ts\frac{\S+1}{2}\right\rc 
     \S < \left\lc\ts\frac{\S+1}{2}\right\rc \S. 
\]

Therefore Eq.(\ref{compspleqn}) translates monochromatic solutions to Eq.(\ref{spleqn}) corresponding to $c \in \big(\S-1,\lc\frac{\S+1}{2}\rc \S\big)$ to monochromatic solutions to Eq.(\ref{spleqn}) corresponding to $c \in \big[\lc\frac{\S+1}{2}\rc \S,\S(\S-1)\big)$.   
\end{Pf}
\vskip 20pt




\begin{thebibliography}{99}

\bibitem{ABEMRS16}
S. D. Adhikari, L. Boza, S. Eliahou, J. M. Mar\'{i}n, M. P. Revuelta and M. I. Sanz, On the $n$-color Rado number for the equation $x_1+x_2+\cdots+x_k+c=x_{k+1}$, {\it Math. Comp.\/} {\bf 85} (2016), 2047--2064. 

\bibitem{ABEMRS17}
S. D. Adhikari, L. Boza, S. Eliahou, J. M. Mar\'{i}n, M. P. Revuelta and M. I. Sanz, On the finiteness of some $n$-color Rado numbers, {\it Discrete Math. \/} {\bf 340} (2017), 39--45. 

\bibitem{BB82}
A. Beutelspacher and W. Brestovansky, Generalized Schur Numbers, {\it Lecture Notes in Mathematics\/} Vol. 969, pp. 30--38. Springer-Verlag, Berlin, 1982. 

\bibitem{BL}
S. A. Burr and S. Loo, On Rado numbers I, unpublished. 

\bibitem{Bro61}
S. L. Brown Jr., Notes on Complete Sequences of Integers, {\it Amer. Math. Monthly\/} {\bf 68} (1961), 557--560.

\bibitem{GS08}
S. Guo and Z.-W. Sun, Determination of two-color Rado number for $a_1x_1+\cdots+a_mx_m=x_0$, {\it J. Combin. Theory Ser. A\/} {\bf 115} (2008), 345--353.  

\bibitem{GTT15}
S. Gupta, J. Thulasi Rangan, and A. Tripathi, The two-colour Rado number for the equation $ax+by=(a+b)z$, {\it Ann. Comb.\/} {\bf 19} (2015), 269--291.

\bibitem{HM97}
H. Harborth and S. Maasberg, Rado numbers for $a(x+y)=bz$, {\it J. Combin. Theory Ser. A\/} {\bf 80} (1997), 356--363. 

\bibitem{HM99}
H. Harborth and S. Maasberg, All two-color Rado numbers for $a(x+y)=bz$, {\it Discrete Math.} {\bf 197/198} (1999), 397--407.

\bibitem{HS05}
B. Hopkins and D. Schaal, All Rado numbers for $\sum_{i=1}^{m-1} a_ix_i=x_m$, {\it Adv. Appl. Math.\/} {\bf 35} (2005), 433--441. 

\bibitem{JS01}
S. Jones and D. Schaal, Some $2$-color Rado numbers, {\it Congr. Numer.} {\bf 152} (2001), 197--199.

\bibitem{JS04}
S. Jones and D. Schaal, Two color Rado number for $x+y+c=kz$, {\it Discrete Math.\/} {\bf 289} (2004), 63--69.

\bibitem{KS01}
W. Kosek and D. Schaal, Rado numbers for the equation $\sum_{i=1}^{m-1} x_i+c=x_m$ for negative values of $c$, {\it Adv. Appl. Math.\/} {\bf 27 (4)} (2001), 805--815.

\bibitem{KSW09}
A. E. K\'{e}zdy, H. S. Snevily, and S. C. White, Generalized Schur Numbers for $x_1+x_2+c=3x_3$, {\it Electron. J. Combin.\/} {\bf 16} (2009), \#105, 13 pp. 

\bibitem{LR14}
B. M. Landman and A. Robertson, Ramsey Theory on the Integers, Second Edition, Student Mathematics Library, {\bf 73}, {\it American Mathematical Society\/}, Providence, RI, 2014. 

\bibitem{Rad33}
R. Rado, Studien zur Kombinatorik, {\it Math. Z.\/} {\bf 36} (1933), 424--480. 

\bibitem{Rad33a}
R. Rado, Verallgemeinerung eines Satzes von van der Waerden mit Anwendungen auf ein Problem der Zahlentheorie, Sonderausg. {\it Sitzungsber. Preuss. Akad. Wiss. Phys. Math. Kl.} {\bf 17} (1933), 1--10.

\bibitem{Rad33b}
R. Rado, Studien zur Kombinatorik, {\it Math. Z.} {\bf 36} (1933), 242--270.

\bibitem{Rad45}
R. Rado, Note on combinatorial analysis, {\it Proc. London Math. Soc.} {\bf 48} (1945), 122--160.

\bibitem{RM08}
A. Robertson and K. Myers, Some two color, four variable Rado numbers, {\it Adv. Appl. Math.\/} {\bf 41} (2008), 214--226.

\bibitem{Sar1}
D. Saracino, The $2$-color Rado number for $x_1+x_2+\cdots+x_{m-1}=ax_m$, {\it arXiv 1207.0432v2\/}

\bibitem{Sar2}
D. Saracino, The $2$-color Rado number for $x_1+x_2+\cdots+x_{m-1}=ax_m$, II, {\it arXiv 1306.0775v1\/}, 4 June 2013, 11 pp. 

\bibitem{Sch93}
D. Schaal, On generalized Schur numbers, {\it Congr. Numer.\/} {\bf 98} (1993), 178--187. 

\bibitem{SZ10}
D. Schaal and M. Zinter, Real number Rado numbers for $\sum_{i=1}^{m-1} x_i +c=x_m$, {\it Congr. Numer.\/} {\bf 202} (2010), 113--118.


\bibitem{Sch16}
I. Schur, \"{U}ber die Kongruenz $x^m+y^m \equiv z^m\mod{p}$, {\it Jahresber. Deutsch. Math. Verein.\/} {\bf 25} (1916), 114--117.  

\end{thebibliography}
\end{document}